\documentclass[conference]{IEEEtran}
\newtheorem{theorem}{Theorem}
\newtheorem{lemma}{Lemma}
\newcommand{\e}{\varepsilon}
\newcommand\numberthis{\addtocounter{equation}{1}\tag{\theequation}}

\usepackage{cite}

\usepackage{amsmath}
\usepackage{amssymb}

\usepackage{algorithm}
\usepackage{algpseudocode}

\begin{document}
%
\title{Adaptive Mirror Descent \\ for Constrained Optimization}

\author{\IEEEauthorblockN{Anastasia Bayandina}
\IEEEauthorblockA{Moscow Institute of Physics and Technology\\
Moscow, Russia\\
Email: anast.bayandina@gmail.com}}

\maketitle

\begin{abstract}
This paper seeks to address how to solve non-smooth convex and strongly convex optimization problems with functional constraints. The introduced Mirror Descent (MD) method with adaptive stepsizes is shown to have a better convergence rate than MD with fixed stepsizes due to the improved constant. For certain types of constraints, the method is proved to generate dual solution. For the strongly convex case, the 'restart' technique is applied.
\end{abstract}

\IEEEpeerreviewmaketitle


\section{Introduction}

Optimizing non-smooth functions with constraints is attracting widespread interest in large-scale optimization and its applications \cite{bib_Shpirko}, \cite{bib_ttd}. There are various methods of solving this kind of optimization problems. The examples of these methods are: bundle-level method \cite{bib_Nesterov}, penalty method \cite{bib_Vasilyev}, \cite{bib_Lan}, Lagrange multipliers method \cite{bib_Boyd}. Among them, Mirror Descent (MD) \cite{bib_Yudin}, \cite{bib_Beck} is viewed as a simple method for non-smooth convex optimization.

In this paper, it is proposed to modify MD so that the stepsizes along with the rate of convergence are no more dependent on the global Lipschitz constant \cite{bib_Juditsky}, but rather on the sizes of the gradients in current points. These sizes are averaged in some sense and substitute the Lipschitz constant. If the constraints can be represented as the maximum of convex functions, which often arises in applications with maximum of many scalar constraints, it is possible to build up the dual solution using the proposed method. The idea of restarts \cite{bib_JuNe} is adopted to construct the algorithm in the case of strongly convex objective and constraints. Both proposed methods are optimal in terms of the lower bounds \cite{bib_Yudin}.

The paper is organized as follows: in Section II we state the problem and notation; in Section III we describe the MD algorithm with adaptive stepsizes and prove the convergence theorem for it; Section IV is focused on the strongly convex case with restarting MD algorithm and theoretical estimates of its convergence; finally, Section V is about duality of the proposed MD method.


\section{Preliminaries and Problem Statement}

Let $E$ be the $n$-dimensional vector space. Let $\lVert\cdot\rVert$ be an arbitrary norm in $E$ and $\lVert\cdot\rVert_*$ be the conjugate norm in $E^*$:
$$\lVert\xi\rVert_* = \max\limits_{x} \big\{ \langle \xi, x \rangle, \lVert x \rVert \leq 1 \big\}.$$

Let $\mathcal{X} \subset E$ be a closed convex set. We consider the two convex functions $f : \mathcal{X} \rightarrow \mathbb{R}$ and $g : \mathcal{X} \rightarrow \mathbb{R}$ to be subdifferentiable and Lipschitz continuous, i.e.
\begin{gather*}
\forall x, y \in \mathcal{X} \hspace{0.2cm} \exists \nabla f(x) : \\
f(y) \geq f(x) + \langle \nabla f(x), y-x \rangle, \hspace{0.3cm}\lVert \nabla f(x) \rVert_* < \infty
\end{gather*}
and the same goes for $g$.

We focus on the problem expressed in the form
\begin{equation}
\label{problem_statement_f}
    f(x) \rightarrow \min\limits_{x\in \mathcal{X}},
\end{equation}
\begin{equation}
\label{problem_statement_g}
    \mathrm{s.t.} \hspace{0.3cm} g(x) \leq 0.
\end{equation}
Denote $x_*$ to be the genuine solution of the problem~(\ref{problem_statement_f}),~(\ref{problem_statement_g}).

Assume that we are equipped with the first-order oracle, which given the point $x\in \mathcal{X}$ returns the values of $\nabla f(x), \nabla g(x),$ and $g(x)$.

Consider $d : \mathcal{X} \rightarrow \mathbb{R}$ to be a distance generating function (d.g.f) which is continuously differentiable and strongly convex, modulus 1, w.r.t. the norm $\lVert\cdot\rVert$, i.e.
$$\forall x, y, \in \mathcal{X} \hspace{0.2cm} \langle d'(x) - d'(y), x-y \rangle \geq \lVert x-y \rVert^2,$$
and assume that
$$\min\limits_{x\in \mathcal{X}} d(x) = d(0).$$
Suppose we are given a constant $\Theta_0$ such that
$$d(x_{*}) \leq \Theta_0^2.$$
Note that if there is a set of optimal points $\mathcal{X}_*$, than we may assume that
$$\min\limits_{x_* \in \mathcal{X}_*} d(x_*) \leq \Theta_0^2.$$
For all $x, y\in \mathcal{X}$ consider the corresponding Bregman divergence
$$V(x, y) = d(y) - d(x) - \langle d'(x), y-x \rangle.$$

For all $x\in \mathcal{X}$, $y\in E^*$ define the proximal mapping operator
$$\mathrm{Mirr}_x (y) = \arg\min\limits_{u\in \mathcal{X}} \big\{ \langle y, u \rangle + V(x, u) \big\}.$$
We make the simplicity assumption, which means that $\mathrm{Mirr}_x (y)$ is easily computable.


\section{Mirror Descent for Constrained Optimization}

The following algorithm is proposed to solve the problem~(\ref{problem_statement_f}),~(\ref{problem_statement_g}).

\begin{algorithm}[H]
\caption{Constraint Mirror Descent}\label{MD}
\begin{algorithmic}[1]
    \Procedure{MD}{$\e, \Theta_0^2, d(\cdot), \mathcal{X}$}
        \State $x^0 \gets \arg\min\limits_{x\in \mathcal{X}} d(x)$
        \State initialize the empty set $I$
        \State $i \gets 0$
        \Repeat
            \If{$g(x^i) \leq \e$}
                \State $M_i \gets \lVert \nabla f(x^i) \rVert_{*}$
                \State $h_i \gets \frac{\e}{M_i^2}$
                \State $x^{i+1} \gets \mathrm{Mirr}_{x^i}(h_i \nabla f(x^i))$
                \State add $i$ to $I$
            \Else
                \State $M_i \gets \lVert \nabla g(x^i) \rVert_{*}$
                \State $h_i \gets \frac{\e}{M_i^2}$
                \State $x^{i+1} \gets \mathrm{Mirr}_{x^i}(h_i \nabla g(x^i))$
            \EndIf
            \State $i \gets i + 1$
        \Until{$\sum\limits_{j \leq i} \frac{1}{M_j^2} \geq \frac{2\Theta_0^2}{\e^2}$}
        \State $\bar{x}^N \gets \frac{\sum\limits_{i\in I} h_i x^i}{\sum\limits_{i\in I} h_i}$
        \State\Return $\bar{x}^N$
    \EndProcedure
\end{algorithmic}
\end{algorithm}

Denote $[N] = \{ i\in \overline{0, N-1} \}$, $J = [N] / I$.

We are going to adopt the following lemma \cite{bib_Nemirovski}.
\begin{lemma}
    Let $f$ be some convex subdifferentiable function over the convex set $\mathcal{X}$. Let the sequence $\{x^i\}$ be defined by the update
    $$ x^{i+1} = \mathrm{Mirr}_{x^i}(h_i \nabla f(x^i)). $$
    Then, for any $x \in \mathcal{X}$
    \begin{equation}
        \begin{gathered}
        h_i \big( f(x^i) - f(x) \big) \leq \\
        \frac{h_i^2}{2} \lVert \nabla f (x^i) \rVert^2_* + V(x^i, x) - V(x^{i+1}, x).
        \end{gathered} \label{Duchi}
    \end{equation}
\end{lemma}

\begin{theorem}\label{th_1}
    The point $\bar{x}^N$ supplied by Algorithm \ref{MD} satisfies
    \begin{equation}
    \label{xxx}
        f(\bar{x}^N) - f(x_{*}) \leq \e, \hspace{0.3cm} g(\bar{x}^N) \leq \e
    \end{equation}
    for the number of oracle calls equal to
    \begin{equation}
        N = \Big\lceil\frac{2M^2 \Theta_0^2}{\e^2}\Big\rceil,
    \end{equation}
    where $M$ is found from
    \begin{equation}
        \frac{N}{M^2} = \sum\limits_{i\in [N]} \frac{1}{M_i^2}.
    \end{equation}
\end{theorem}

\begin{IEEEproof}
    By the definition of $\bar{x}^N$ and the convexity of $f$,
    \begin{equation}\label{th1_eq1}
        \sum\limits_{i\in I} h_i f(\bar{x}^N) - f(x_{*}) \leq \sum\limits_{i\in I} h_i \big( f(x^i) - f(x_{*}) \big).
    \end{equation}
    
    Using (\ref{Duchi}) and the definitions of the stepsizes, consider the summation 
    \begin{gather*}
        \sum\limits_{i\in I} h_i \big( f(x^i) - f(x_{*}) \big) + \sum\limits_{i\in J} h_i \big( g(x^i) - g(x_{*}) \big) \leq \\
        \sum\limits_{i\in I} \frac{h_i^2 M_i^2}{2} + \sum\limits_{i\in J} \frac{h_i^2 M_i^2}{2} + \\
        \sum\limits_{i\in [N]} \big( V(x^i, x_{*}) - V(x^{i+1}, x_{*}) \big) \leq \frac{\e}{2} \sum\limits_{i\in [N]} h_i + \Theta_0^2. \numberthis \label{th1_eq2}
    \end{gather*}
    
    Since for $i\in J$
    $$g(x^i) - g(x_{*}) \geq g(x^i) > \e,$$
    recalling (\ref{th1_eq1}), we get
    \begin{gather*}
        \sum\limits_{i\in I} h_i \big( f(\bar{x}^N) - f(x_{*}) \big) < \frac{\e}{2} \sum\limits_{i\in [N]} h_i - \e\sum\limits_{i\in J} h_i + \Theta_0^2 = \\
        \e\sum\limits_{i\in I} h_i - \frac{\e^2}{2} \sum\limits_{i\in [N]} \frac{1}{M_i^2} + \Theta_0^2 \leq \e\sum\limits_{i\in I} h_i.
    \end{gather*}
    As long as the inequality is strict, the case of the empty $I$ is impossible.
    
    For $i\in I$ holds $g(x^i) \leq \e$. Then, by the definition of $\bar{x}^N$ and the convexity of $g$,
    \begin{gather*}
        \sum\limits_{i\in I} h_i g(\bar{x}^N) \leq \sum\limits_{i\in I} h_i g(x^i) \leq \e\sum\limits_{i\in I} h_i.
    \end{gather*}    
\end{IEEEproof}

It is worth mentioning that the constant $M$ is somewhat the average of all subgradient norms in particular points instead of being the Lipschitz constant biggest possible over the set~$\mathcal{X}$.


\section{Restarting Mirror Descent}

In this section we assume that $f$ and $g$ in the problem~(\ref{problem_statement_f}),~(\ref{problem_statement_g}) are $\mu$-strongly convex on $\mathcal{X}$, i.e.
$$\forall x, y, \in \mathcal{X} \hspace{0.2cm} f(y) \geq f(x) + \langle \nabla f(x), y-x \rangle + \frac{\mu}{2} \lVert y-x \rVert^2,$$
and the same goes for $g$.

Also the d.g.f is assumed to be bounded on the unit ball, that is,
\begin{equation}\label{d}
\forall x\in \mathcal{X} : \lVert x \rVert \leq 1 \hspace{0.2cm} d(x) \leq \frac{\omega_n}{2},
\end{equation}
where $\omega_n$ is some dimension-dependent constant which in most setups asymptotically behaves as $O\big(\log(n)\big)$ \cite{bib_Nemirovski}.

Suppose we are given a constant $R_0$ such that
\begin{equation}\label{R_0}
\lVert x_0 - x_* \rVert^2 \leq R_0^2, \text{ where } x_0 = \arg\min\limits_{x\in \mathcal{X}} d(x).
\end{equation}

The following algorithm is proposed to solve the problem~(\ref{problem_statement_f}),~(\ref{problem_statement_g}) in the case of strong convexity \cite{bib_JuNe}.

\begin{algorithm}[H]
\caption{Restarting Mirror Descent}\label{Restart_MD}
\begin{algorithmic}[1]
    \Procedure{RestartMD}{$\e, \omega_n, R_0, \mathcal{X}$}
        \State $x_0 \gets \arg\min\limits_{x\in \mathcal{X}} d(x)$
        \State $d_0(x) \gets d(\frac{x - x_0}{R_0})$
        \State $K \gets \log_2 \frac{\mu R_0^2}{2\e}$
        \For{$k \gets 1, K$}
            \State $R_k^2 \gets R_0^2 \cdot 2^{-k}$
            \State $\e_k \gets \frac{\mu R_k^2}{2}$
            \State $x_k \gets \mathrm{MD}(\e_k, \frac{\omega_n R_k^2}{2}, d_{k-1}(\cdot), \mathcal{X})$
            \State $d_k(x) \gets d(\frac{x - x_k}{R_k})$
        \EndFor
        \State\Return $x_K$
    \EndProcedure
\end{algorithmic}
\end{algorithm}

At each iteration $k$ of the loop the algorithm performs the restart: it calls the $\mathrm{MD}$ procedure described in the previous section with some accuracy $\e_k$ which becomes smaller for each next restart.

Denote by $N_1, \dots, N_K$ the numbers of oracle calls at each restart in Algorithm \ref{Restart_MD} and by $[N_1], \dots, [N_K]$ the corresponding sets of indices.

Further for the sake of brevity we accept the following statement without proof.

\begin{lemma}\label{lemma2}
    Suppose $f$ and $g$ are $\mu$-strongly convex functions w.r.t. the norm $\lVert\cdot\rVert$ and $x_*$ is the genuine solution of the problem~(\ref{problem_statement_f}),~(\ref{problem_statement_g}). Then if for some $x\in \mathcal{X}$
    $$f(x) - f(x_*) \leq \e, \hspace{0.3cm} g(x) \leq\e,$$
    then
    $$\frac{\mu}{2} \lVert x - x_* \rVert \leq \e.$$
\end{lemma}

Now we are ready to prove the following

\begin{theorem}
    The point $x_K$ supplied by Algorithm \ref{Restart_MD} satisfies
    \begin{equation}
        f(x_K) - f(x_{*}) \leq \e, \hspace{0.3cm} g(x_K) \leq \e
    \end{equation}
    for the total number of oracle calls equal to
    \begin{equation}
        N = N_1 + \dots + N_K = \Big\lceil \frac{4M^2 \omega_n}{\mu\e} \Big\rceil,
    \end{equation}
    where $M$ is found from
    \begin{equation}
        \frac{N}{M^2} = \sum\limits_{i\in [N_1]} \frac{1}{M_i^2} + \dots + \sum\limits_{i\in [N_K]} \frac{1}{M_i^2}.
    \end{equation}
\end{theorem}
\begin{IEEEproof}
    Observe \cite{bib_Juditsky} that the function $d_{k-1}(x)$ defined in Algorithm \ref{Restart_MD} is 1-strongly convex w.r.t. the norm $\lVert\cdot\rVert / R_{k-1}$. The conjugate of this norm is $R_{k-1} \lVert\cdot\rVert_*$. It means that at each restart the actual Lipschitz constants are $M_i R_{k-1}$. Then, by (\ref{d}) and (\ref{R_0}) at the end of the first restart we obtain
    $$\sum\limits_{i\in [N_1]} \frac{1}{M_i^2 R_0^2} =  \frac{\omega_n}{\e_1^2} \geq \frac{2d_0(x_*)}{\e_1^2},$$
    which by Theorem \ref{th_1} guarantees the $\e_1$-solution of the problem.
    
    Further, by Lemma \ref{lemma2}, after the $(k-1)$th restart it holds that
    $$\lVert x_{k-1} - x_* \rVert^2 \leq \frac{2\e_{k-1}}{\mu} = R_{k-1}^2.$$
    Due to the choice of the d.g.f. $d_{k-1}(x)$, the starting point of the $k$th restart is $x_{k-1}$ and
    
    $$\sum\limits_{i\in [N_k]} \frac{1}{M_i^2 R_{k-1}^2} =  \frac{\omega_n}{\e_k^2} \geq \frac{2d_{k-1}(x_*)}{\e_k^2}.$$
    In that way we have justified the redefinition of the d.g.f. and the 'distance' argument of the $\mathrm{MD}$ procedure.

    After the $k$th restart by the definition of $R_k$ and $\e_k$ we obtain
    $$\sum\limits_{i\in [N_k]} \frac{1}{M_i^2} = \frac{\omega_n R_k^2}{\e_k^2} = \frac{4\omega_n}{\mu^2 R_0^2} \cdot 2^k.$$
    Thus, for the whole $\mathrm{RestartMD}$ procedure considering the definition of $K$
    \begin{gather*}
        \sum\limits_{i\in [N_1]} \frac{1}{M_i^2} + \dots + \sum\limits_{i\in [N_K]} \frac{1}{M_i^2} = \frac{4\omega_n}{\mu^2 R_0^2} \cdot (2 + \dots + 2^K) = \\
        \frac{8\omega_n}{\mu^2 R_0^2} \cdot 2^K - \frac{8\omega_n}{\mu^2 R_0^2} \leq \frac{4\omega_n}{\mu\e}.
    \end{gather*}
\end{IEEEproof}

Note that due to Lemma \ref{lemma2} the argument $x_k$ converges to $x_*$ along with the function, which is a typical property of strongly convex optimization.


\section{Dual Problem Solution}

Following \cite{bib_arxiv}, in this section we regard the problem of the type (\ref{problem_statement_f}),~(\ref{problem_statement_g}) where the constraints appear in the form
\begin{equation}\label{max_constraints}
g(x) = \max\limits_{m \in \overline{1, \mathcal{M}}} \big\{ g_m(x) \big\}.
\end{equation}
Consider the dual problem
\begin{equation}\label{dual}
    \varphi (\lambda) = \min\limits_{x\in \mathcal{X}} \Big\{ f(x) + \sum\limits_{m=1}^{\mathcal{M}} \lambda_m g_m(x) \Big\} \rightarrow \max\limits_{\lambda \geq 0}.
\end{equation}

Denote $\lambda_* = (\lambda_{*1}, \dots, \lambda_{*\mathcal{M}})$ to be the genuine solution of (\ref{dual}). Then, by the weak duality property \cite{bib_Boyd} we have
$$\Delta (x_*, \lambda_*) \stackrel{\text{def}}{=} f(x_*) - \varphi(\lambda_*) \geq 0.$$
Assume that Slater's condition holds, i.e. there exists $x\in \mathcal{X}$ such that $g(x) < 0$. This ensures strong duality $\Delta (x_*, \lambda_*) = 0$. It means that if the algorithm is able to generate the dual solution $\bar{\lambda}^N = (\bar{\lambda}_1^N, \dots, \bar{\lambda}_{\mathcal{M}}^N)$ of the problem (\ref{problem_statement_f}),~(\ref{problem_statement_g}) with (\ref{max_constraints}), the accuracy of this solution can be estimated via the size of the duality gap $\Delta (\bar{x}^N, \bar{\lambda}^N)$.

As long as the constraints are of the form (\ref{max_constraints}), we can define the function
\begin{equation}\label{m}
    m(i) : [N] \rightarrow \{ 1, \dots, \mathcal{M} \}, \hspace{0.3cm} g(x^i) = g_{m(i)}(x^i).
\end{equation}

\begin{theorem}
    Consider Algorithm \ref{MD} and define dual Lagrange multipliers as
    \begin{equation}\label{lam}
        \bar{\lambda}_m^N = \frac{1}{\sum\limits_{i\in I} h_i} \sum\limits_{i\in J} h_i \mathbb{I}\big[ m(i) = m \big],
    \end{equation}
    where
    $$\mathbb{I}[x] = \begin{cases}
        1, \text{ if $x$ is true}, \\
        0, \text{ otherwise.}
    \end{cases}$$
    Then, the point $\bar{x}^N$ supplied by Algorithm \ref{MD} satisfies
    \begin{equation}
        \Delta(\bar{x}^N, \bar{\lambda}^N) = f(\bar{x}^N) - \varphi(\bar{\lambda}^N) \leq \e, \hspace{0.3cm} g(\bar{x}^N) \leq \e
    \end{equation}
    for the number of oracle calls equal to
    \begin{equation}
        N = \Big\lceil\frac{2M^2 \Theta_0^2}{\e^2}\Big\rceil,
    \end{equation}
    where $M$ is found from
    \begin{equation}
        \frac{N}{M^2} = \sum\limits_{i\in [N]} \frac{1}{M_i^2}.
    \end{equation}
\end{theorem}
\begin{IEEEproof}
    Combining (\ref{th1_eq1}) and (\ref{th1_eq2}) together with (\ref{m}) we obtain
    \begin{gather*}
        \sum\limits_{i\in I} h_i \big( f(\bar{x}^N) - f(x_*) \big) \leq \frac{\e}{2} \sum\limits_{i\in [N]} h_i + \Theta_0^2 - \\
        \sum\limits_{i\in J} h_i \big( g(x^i) - g(x_*) \big) = \frac{\e}{2} \sum\limits_{i\in [N]} h_i + \Theta_0^2 - \\
        \sum\limits_{i\in J} h_i \big( g_{m(i)}(x^i) - g_{m(i)}(x_*) \big) \leq \\
        \e \sum\limits_{i\in I} h_i + \sum\limits_{i\in J} \sum\limits_{m=1}^{\mathcal{M}} h_i \mathbb{I}\big[ m(i) = m \big] g_m(x_*).
    \end{gather*}
    Recalling (\ref{lam}) and rearranging the terms,
    \begin{gather*}
        \sum\limits_{i\in I} h_i f(\bar{x}^N) \leq \e\sum\limits_{i\in I} h_i + \\
        \sum\limits_{i\in I} h_i \Big( f(x_*) + \sum\limits_{m=1}^{\mathcal{M}} \bar{\lambda}_m^N g_m(x_*) \Big) = \\
        \e\sum\limits_{i\in I} h_i + \sum\limits_{i\in I} h_i \min\limits_{x\in \mathcal{X}} \Big\{ f(x) + \sum\limits_{m=1}^{\mathcal{M}} \bar{\lambda}_m^N g_m(x) \Big\} = \\
        \e\sum\limits_{i\in I} h_i + \sum\limits_{i\in I} h_i \varphi (\bar{\lambda}^N).
    \end{gather*}
\end{IEEEproof}


\section{Conclusion}
We proved MD algorithm with adaptive stepsizes to achieve optimal rates in both convex and strongly convex cases with the improved Lipschitz constant. For the problems with constraints in the form of maximum of convex functions we showed the duality of the method. However, it still remains open whether it is possible to construct high probability bounds for adaptive steps in the case of stochastic oracle.


\section*{Acknowledgment}

The author gratefully acknowledges the help and valuable discussion kindly provided by Dr. Gasnikov.

This research was funded by Russian Science Foundation (project 17-11-01027).


\end{document}